\documentclass[12pt,draft]{amsart}
\usepackage{amsmath,amsthm,latexsym,amscd,amsbsy,amssymb,fleqn,leqno}
\chardef\bslash=`\\ % p. 424, TeXbook

\makeatletter
\def\verbatim{\interlinepenalty\@M \@verbatim
  \leftskip\@totalleftmargin\advance\leftskip2pc
  \frenchspacing\@vobeyspaces \@xverbatim}
\makeatother
\makeatletter
  \def\dgt@k{\dg@DX=-3 \dg@DY=2 \dg@SIZE=3} 
\makeatother

\makeatletter
  \def\dgt@kk{\dg@DX=3 \dg@DY=-1 \dg@SIZE=3}%
\makeatother

\theoremstyle{plain}
\newtheorem{thm}{Theorem}[section]
\newtheorem{cor}[thm]{Corollary}
\newtheorem{lem}[thm]{Lemma}
\newtheorem{pro}[thm]{Proposition}

\theoremstyle{definition}

\numberwithin{equation}{section}

\newcounter{rmnum}

\def\symbolnote#1#2{\let\thefootn=\thefootnote%
\renewcommand{\thefootnote}{\fnsymbol{footnote}}%
\footnotemark[#1]%
\footnotetext[#1]{#2}%
\let\thefootnote=\thefootn
}

\newfont{\bbb}{msbm10 scaled \magstep1}
\newfont{\bbc}{msbm8 scaled \magstep0}

\newcommand{\N}{\mbox{\bbb N}}

\newcommand{\uin}{\mbox{\bbb I}}

\begin{document}

%%%%%%% Begin Topmatter %%%%%%%%%%

\title[On finite-to-one maps]{On finite-to-one maps}
\author{H. Murat Tuncali}
\address{Department of Mathematics,
Nipissing University,
100 College Drive, P.O. Box 5002, North Bay, ON, P1B 8L7, Canada}
\email{muratt@nipissingu.ca}
\thanks{The first author was partially supported by NSERC grant.}

\author{Vesko Valov}
\address{Department of Mathematics, Nipissing University,
100 College Drive, P.O. Box 5200, North Bay, ON, P1B 8L7, Canada}
\email{veskov@nipissingu.ca}
\thanks{The second author was partially supported by Nipissing University Research Council Grant.}

\keywords{finite-to-one maps, dimension, set-valued maps} 
\subjclass{Primary: 54F45; Secondary: 55M10, 54C65.}
 
%%%%%%% End topmatter %%%%%%%%%

\begin{abstract}{Let $f\colon X\to Y$ be a $\sigma$-perfect $k$-dimensional surjective map of metrizable spaces such that $\dim Y\leq m$. It is shown that,  for every integer $p\geq 1$ there exists a dense $G_{\delta}$-subset 
${\mathcal H}(k,m,p)$ of $C(X,\uin^{k+p})$ with the source limitation topology such that   each fiber of 
$f\triangle g$, $g\in{\mathcal H}(k,m,p)$, contains at most $\max\{k+m-p+2,1\}$ points.  This result provides a proof of Hypothesis 1 and Hypothesis 2 from \cite{bfm:00}}.   
\end{abstract}

\maketitle
\markboth{H. M.~Tuncali and V.~Valov}{On finite-to-one maps}

%%%%%%%%%%%%%%%%%%%%%%%%%%%%%%%%%%%%%%%%%%%%%%%%%%%%%%%%%%%%

\section{Introduction}

This paper is inspired by the following hypotheses of S. Bogatyi, V. Fedorchuk and J. van Mill  \cite{bfm:00}. 

Let $f\colon X\to Y$ be at most a $k$-dimensional map between compact metric spaces with $\dim Y\leq m$. Then: (1) there exists a map $h\colon X\to\uin^{m+2k}$ such that $f\triangle h\colon X\to Y\times\uin^{m+2k}$ is at most 2-to-one provided $k\geq 1$;
(2) there exists a map  $h\colon X\to\uin^{m+k+1}$ such that $f\triangle h\colon X\to Y\times\uin^{m+k+1}$ is at most $(k+1)$-to-one.

Next theorem provides a solution of these two problems (by a $k$-dimensional map we mean a map with all fibers being at most $k$-dimensional).
 
\begin{thm}
Let $f\colon X\to Y$ be a $\sigma$-perfect $k$-dimensional surjective map of metrizable spaces such that $\dim Y\leq m$.  Then,  for every integer $p\geq 1$, the set ${\mathcal H}(k,m,p)$ consisting of all $g\in C(X,\uin^{k+p})$ such that each fiber of the map
$f\triangle g\colon X\to Y\times\uin^{k+p}$ contains at most $\max\{k+m-p+2, 1\}$ points is dense and
$G_{\delta}$ in $C(X,\uin^{k+p})$ with respect to the source limitation topology.
\end{thm}

Observe that stronger forms of Hypothesis 1 and Hypothesis 2 follow from Theorem 1.1 when $p=m+k$ and $p=m+1$, respectively.  If $p=m+k+1$, then all maps $f\triangle g$, $g\in{\mathcal H}(k,m,p)$ are one-to-one.  So, Theorem 1.1 implies
\cite[Theorem 7.3]{bp:96} and the metrizable case of  \cite[Theorem 1.1(i)]{tv1:02}. 
When both $X$ and $Y$ are compact, $k=0$ and $p=1$ Theorem 1.1 was established by   
 M. Levin and W. Lewis \cite[Proposition 4.4]{ll:02}. This result is one of the ingredients of our proof, another one is a selection theorem proven by V. Gutev and the second author \cite[Theorem 1.2]{gv:01}.

Recall that, $\dim f=\sup\{\dim f^{-1}(y):y\in Y\}$ is the dimension of $f$. We say that 
a surjective map $f\colon X\to Y$ is called $\sigma$-perfect if $X$ is the union of countably many closed sets $X_i$ such that 
each restriction $f|X_i\colon X_i\to f(X_i)$ is a perfect map. By $C(X,M)$ we denote the set of all continuous maps from
$X$ into $M$. If $(M,d)$ is a metric space,  then the source limitation topology on $C(X,M)$ is defined in the following way:
a subset $U\subset C(X,M)$ is open in $C(X,M)$ with respect to
the source limitation topology provided for every $g\in U$ there exists 
a continuous function $\alpha\colon X\to (0,\infty)$ such that $\overline{B}(g,\alpha)\subset U$. Here, $\overline{B}(g,\alpha)$ denotes the set  
$\{h\in C(X,M):d(g(x),h(x))\leq\alpha (x)\hbox{}~~\mbox{for each 
$x\in X$}\}$. 
The source limitation topology doesn't depend on the metric $d$ if $X$ is paracompact \cite{nk:69}  and $C(X,M)$ with this topology has Baire property provided $(M,d)$ is a complete metric space \cite{jm:75}.  Moreover, if $d$ is a bounded metric on $M$ and $X$ is compact, then the source limitation topology coincides with the uniform convergence topology generated by $d$.

The paper is organized as follows. In Section 2 we prove the  special case of Theorem 1.1 when both $X$ and $Y$ are compact.   The final proof is accomplished in Section 3.

All maps are assumed to be continuous and all function spaces, if not explicitely stated otherwise, are equipped with the source limitation topology.  Everywhere in this paper by an  $n$-to-one map, where $n\geq 1$ is an integer, we mean a map  with all fibers containing at most $n$ points.

%%%%%%%%%%%%%%%%%%%%%%%%%%%%%%%%%

%%%%%%%%%%%%%%%%%%%%%%%%%%%%%%
%%%%%%%%%%%%%%%%%%%%%%%%%%%%%%
\vspace{0.5cm}

\section{Proof of Theorem 1.1 - the compact case}

Let  $\omega$ be an open cover of the space $X$, $m\in\N$ and $H\subset X$. 
We say that a map
$g\colon H\to Z$ is an $(m,\omega)$-map if every $z\in g(H)$ has a neighborhood $V_z$ in $Z$ such that $g^{-1}(V_z)$ can be covered by $m$ elements of $\omega$. We also agree to denote by $cov(M)$ the family of all open covers of $M$. 

Suppose $f\colon X\to Y$ is a surjective map, $\omega\in cov(X)$ and $n, m\in\N$. Then,
we denote by $C(X,Y\times\uin^n,f)$
the set of all maps $h\colon X\to Y\times\uin^n$ such that $\pi_Y\circ h=f$, where
$\pi_Y\colon Y\times\uin^p\to Y$ is the projection. For any $K\subset X$, $C_{(m,\omega)}(X|K,Y\times\uin^n,f)$ stands for the set of all $h\in C(X,Y\times\uin^n,f)$
with $h|K$ being an $(m,\omega)$-map  and 
$C_{(m,\omega)}(X|K,\uin^n)$ consists of all $g\in C(X,\uin^n)$ such that 
$f\triangle g\in C_{(m,\omega)}(X|K,Y\times\uin^n,f)$. In case $K=X$ we simply write 
$C_{(m,\omega)}(X,Y\times\uin^n,f)$ (resp., $C_{(m,\omega)}(X,\uin^n)$) instead of
$C_{(m,\omega)}(X|X,Y\times\uin^n,f)$ (resp., $C_{(m,\omega)}(X|X,\uin^n)$).

\begin{pro}
Let $f\colon X\to Y$ be a surjection between metrizable spaces and $\{X_i\}$ a sequence of closed subsets of $X$ such that each restriction $f|X_i$ is a perfect map.  
Then for any positive integers $m$ and $p$ the set 
${\mathcal A}(m,p)=\{g\in C( X,\uin^p): (f\triangle g)|(\cup_{i=1}^{\infty}X_i)\mbox{}~\mbox{is $m$-to-one}\}$ is $G_{\delta}$ in $C(X,\uin^p)$.
\end{pro}

\begin{proof}
We need few lemmas, in all these lemmas we suppose that $X$, $Y$ and $f$ are as in Proposition 2.1 and $\omega\in cov(X)$. 

\begin{lem}
Let $f$ be a perfect map and $g\in C_{(m,\omega)}(X|f^{-1}(y),\uin^p)$ for some $y\in Y$. Then there exists a neighborhood $U_y$ of $y$ in $Y$ such that the restriction $g|f^{-1}(U_y)$ is an $(m,\omega)$-map.
\end{lem}

\begin{proof}
Obviously, $g\in C_{(m,\omega)}(X|f^{-1}(y),\uin^p)$ implies that 
$g|f^{-1}(y)$ is an $(m,\omega)$-map. Hence, for every $x\in f^{-1}(y)$ there exists an open neighborhood $V_{g(x)}$ of $g(x)$ in $\uin^p$ such that $g^{-1}(V_{g(x)})\cap f^{-1}(y)$ can be covered by $m$ elements of $\omega$ whose union is denoted by
 $W_x$.
Therefore, for every $x\in f^{-1}(y)$ we have
$(f\triangle g)^{-1}(f(x),g(x))=f^{-1}(y)\cap g^{-1}(g(x))\subset W_x$ and, since $f\triangle g$ is a closed map, there exists an open neighborhood $H_x=U_{y}^x\times G_x$ of $(y,g(x))$ in $Y\times\uin^p$ with $S_x=(f\triangle g)^{-1}(H_x)\subset W_x$. Next, choose finitely many points $x(i)\in f^{-1}(y)$, $i=1,2,..n$, 
such that $f^{-1}(y)\subset \bigcup_{i=1}^{i=n}S_{x(i)}$. Using that $f$ is a closed map we can find a neighborhood $U_y$ of $y$ in $Y$ such that $U_y\subset\bigcap_{i=1}^{i=n}U_{y}^{x(i)}$ and $f^{-1}(U_y)\subset\bigcup_{i=1}^{i=n}S_{x(i)}$. Let show that $g|f^{-1}(U_y)$ is an $(m,\omega)$-map. Indeed, if $z\in f^{-1}(U_y)$, then $z\in S_{x(j)}$ for some $j$ and $g(z)\in G_{x(j)}$ because $S_{x(j)}=f^{-1}(U_{y}^{x(j)})\cap g^{-1}(G_{x(j)})$. Consequently, 
$f^{-1}(U_{y})\cap g^{-1}(G_{x(j)})\subset S_{x(j)}\subset W_{x(j)}$. Therefore, $G_{x(j)}$ is a neighborhood of $g(z)$ such that $f^{-1}(U_{y})\cap g^{-1}(G_{x(j)})$ is covered by $m$ elements of $\omega$.
\end{proof}

\begin{cor}
If $f$ is perfect and $g\in C_{(m,\omega)}(X|f^{-1}(y),\uin^p)$ for every $y\in Y$, then  
$g\in C_{(m,\omega)}(X,\uin^p)$.
\end{cor}

\begin{proof}
 By Lemma 2.2, for any $x\in X$ there exists
a neighborhood $U_y$ of $y=f(x)$ in $Y$ 
such that $g|f^{-1}(U_y)$ is an $(m,\omega)$-map. 
So, we can find a neighborhood $G_{x}$ of $g(x)$ in $\uin^p$ with  
$f^{-1}(U_{y})\cap g^{-1}(G_{x})$ being covered by $m$ elements of $\omega$.  But $f^{-1}(U_{y})\cap g^{-1}(G_{x})$ equals to
$(f\triangle g)^{-1}(U_y\times G_x)$. Hence, 
$f\triangle g$ is an $(m,\omega)$-map. 
\end{proof} 

\begin{lem}
For any closed $K\subset X$ the set 
$C_{(m,\omega)}(X|K,\uin^p)$ is open in $C(X,\uin^p)$  provided $f$ is perfect.  
\end{lem}

\begin{proof}
The proof of this lemma follows the same scheme as the proof of \cite[Lemma 2.5]{tv:02}, we apply now Lemma 2.2 instead of Lemma 2.3 from \cite{tv:02}.  
\end{proof}

Let finish the proof of Proposition 2.1. We can suppose that the sequence $\{X_i\}$ is increasing and fix a sequence $\{\omega_i\}\subset cov(X)$ such that $mesh(\omega_i)\leq i^{-1}$ for every $i$. Denote by $\pi_i\colon C(X,\uin^p)\to C(X_i,\uin^p)$, $\pi_i(g)=g|X_i$, the restriction maps.  
By Lemma 2.4, every set 
${\mathcal B}_{ij}$, $i,j\in\N$, consisting of all $h\in C(X_i,\uin^p)$ with $(f|X_i)\triangle h$ being an $(m,\omega_j)$-map is open in $C(X_i,\uin^p)$. So are the sets ${\mathcal A}_{ij}=(\pi_i)^{-1}({\mathcal B}_{ij})$ in $C(X,\uin^p)$  because each $\pi_i$ is continuous.  It is easily seen that  the intersection of all ${\mathcal A}_{ij}$ is exactly the set  ${\mathcal A}(m,p)$. Hence, ${\mathcal A}(m,p)$ is a $G_{\delta}$-subset of $C(X,\uin^p)$. 
\end{proof}

\begin{cor}
Theorem $1.1$ follows from the validity of its special case when $p\leq k+m+1$
\end{cor}

\begin{proof}
Suppose $p=m+k+1+n$, where $n\geq 1$.  Then $C(X,\uin^{k+p})$ is homeomorphic to the product 
$C(X, \uin^{2k+m+1})\times C(X,\uin^n)$ and let $\pi\colon C(X,\uin^{k+p})\to C(X,\uin^{2k+m+1})$ denote the projection.  According to our assumption, the set $\mathcal A=\{h\in C(X,\uin^{2k+m+1}): f\triangle h \mbox{}~\mbox{is one-to-one}\}$ is dense in $C(X,\uin^{2k+m+1})$, so is the set  $\pi^{-1}(\mathcal A)$ in $C(X,\uin^{k+p})$. Since $\max\{k+m-p+2,1\}=1$, 
${\mathcal H}(k,m,p)$ consists of one-to-one maps. Hence $\pi^{-1}(\mathcal A)\subset {\mathcal H}(k,m,p)$. The last inclusion 
yields that ${\mathcal H}(k,m,p)$ is dense in $C(X,\uin^{k+p})$. It only remains to observe that, by Proposition 2.1, ${\mathcal H}(k,m,p)$ is $G_{\delta}$ in $C(X,\uin^{k+p})$. 
\end{proof}

The remaining part of this section is devoted to the proof of next proposition which, in combination with Proposition 2.1 and Corollary 2.5  provides a proof of Theorem 1.1 when both $X$ and $Y$ are compact.

\begin{pro}
The set ${\mathcal H}(k,m,p)$ in Teorem $1.1$ is dense in $C(X,\uin^{k+p})$ provided both $X$ and $Y$ are compact metric spaces and $p\leq m+k+1$.
\end{pro} 

\begin{proof}
Let first show that the proof of this proposition can be reduced to the proof of its special case when $k=0$. Indeed, suppose Proposition 2.5 is valid for $k=0$ and every positive $p$ with $p\leq m+1$.
Fix $\epsilon>0$ and $h\in C(X,\uin^{k+p})$, where $k\geq 0$ and $1\leq p\leq m+k+1$. Then $h=h_1\triangle h_2$ with $h_1\in C(X,\uin^k)$ and $h_2\in C(X,\uin^p)$. By \cite{bp:96}, there exists $g_1\in C(X,\uin^k)$ such that $f\triangle g_1\colon X\to Y\times\uin^k$ is a 0-dimensional map and $g_1$ is $\displaystyle\frac{\epsilon}{2}$-close to $h_1$. Then, applying our assumption to the map $f\triangle g_1$, we can find $g_2\in C(X,\uin^p)$ which is 
$\displaystyle\frac{\epsilon}{2}$-close to $h_2$ and such that $(f\triangle g_1)\triangle g_2$ is a
$(k+m-p+2)$-to-one map. It remains only to observe that the map $g=g_1\triangle g_2\in C(X,\uin^{k+p})$ is $\epsilon$-close to $h$ and $f\triangle g$ is a $(k+m-p+2)$-to-one map.

So, the following statement, which is denoted by $\Sigma(m,p)$, will complete the proof:\\ {\em Let $f\colon X\to Y$ be a 0-dimensional surjection between compact metrizable spaces with $\dim Y\leq m$. Then , for every positive integer $p\leq m+1$ the set  
${\mathcal H}(0,m,p)=\{g\in C( X,\uin^p): f\triangle g\mbox{}~\mbox{is $(m-p+2)$-to-one}\}$ is dense in $C(X,\uin^p)$}.\\
We are going to prove $\Sigma(m,p)$ by induction with respect to $p$. The statement $\Sigma(m,1)$ was proved by M. Levin and  W. Lewis \cite[Proposition 4.4]{ll:02}. Assume that $\Sigma(m,p)$ holds for any $p\leq n$ and $m\geq p-1$, where $n\geq 1$, and let prove the validity of $\Sigma(m,n+1)$.    
We need to show that for fixed $m$ with $n\leq m$, $h^*\in C(X,\uin^{n+1})$ and $\epsilon>0$ there exists $g^*\in{\mathcal H}(0,m,n+1)$ which is $\epsilon$-close to $h^*$. 
To this end, we represent $h^*$ as $h^*_1\triangle h^*_2$, where $h^*_1\in C(X,\uin^n)$ and $h^*_2\in C(X,\uin)$. Next,  
 we use an idea from the proof of \cite[Theorem 5]{bfm:00}. By Urysohn's  decomposition theorem (see \cite[Theorem 1.5.7]{re:95}), there exists an $F_{\sigma}$-subset $Y_0\subset Y$ such that $\dim Y_0\leq m-1$ and $\dim (Y\backslash Y_0)=0$. 
Let $Y_0$ be the union of an increasing sequence of closed sets $Y_i\subset Y$, $i\geq 1$,  and $X_i=f^{-1}(Y_i)$, $i\geq 0$.  Obviously, $\dim Y_i\leq m-1$, $i\geq 1$ and $n\leq (m-1)+1$. Thus, according to our inductive hypothesis we can apply $\Sigma(m-1,n)$  for the maps $f_i=f|X_i\colon X_i\to Y_i$, $i\geq 1$, to conclude that 
each set ${\mathcal B}_i=\{g\in C(X_i,\uin^n): f_i\triangle g\mbox{}~\mbox{is $(m-n+1)$-to-one}\}$  is dense in $C(X_i,\uin^n)$.  So are the sets  ${\mathcal A}_i=(\pi_i)^{-1}({\mathcal B}_i)$ in $C(X,\uin^n)$ because the restriction maps $\pi_i\colon C(X,\uin^n)\to C(X_i,\uin^n)$ are open and surjective. 
On the other hand, by Proposition 2.1, each of the sets ${\mathcal A}_i$ is $G_{\delta}$ in $C(X,\uin^n)$.  Hence, the set 
${\mathcal A}_0=\cap_{i=1}^{\infty}{\mathcal A}_i$ is dense and $G_{\delta}$ in $C(X,\uin^n)$ and obviously, it consists of  all maps $q\in C(X,\uin^n)$ such that $(f\triangle q)|X_0$ is $(m-n+1)$-to-one.  Therefore, there exists $g^*_1\in{\mathcal A}_0$ which is $\displaystyle\frac{\epsilon}{2}$-close to $h^*_1$. Consider the map $f\triangle g^*_1\colon X\to Y\times\uin^n$ and the set 
$D=\{z\in Y\times\uin^n: |(f\triangle g^*_1)^{-1}(z)|\geq m-n+2\}$. By \cite[Lemma 4.3.7]{re:95}, $D\subset Y\times\uin^n$ is $F_{\sigma}$.  Then $H=\pi_Y(D)$ doesn't meet $Y_0$ because of the choice of $g^*_1$, where $\pi_Y\colon Y\times\uin^n\to Y$ denotes the projection.  Hence, $H$ is 0-dimensional and $\sigma$-compact, so is the set $K=f^{-1}(H)$ ( 0-dimensionality of $K$ follows by the Hurewicz theorem on dimension-lowering mappings, see \cite[Theorem 1.12.4]{re:95}). Representing $H$ as the union of an increasing sequence of closed sets $H_i\subset Y$  and applying $\Sigma(0,1)$  for any of the maps $f|K_i$ , where $K_i=f^{-1}(H_i)$, we can conclude (as we did for the set ${\mathcal A}_0$ above) that the set ${\mathcal F}$ of all maps $q\in C(X,\uin)$ with $(f\triangle q)|K$ one-to-one  is dense and $G_{\delta}$ in $C(X,\uin)$. Consequently, there exists $g^*_2\in{\mathcal F}$ which is $\displaystyle\frac{\epsilon}{2}$-close to $h^*_2$. Then $g^*=g^*_1\triangle g^*_2$ is $\epsilon$-close to $h^*$.  It follows from the definition of the set $D$ and the choice of the maps $g^*_1, g^*_2$ that $f\triangle g^*$ is $(m-n+1)$-to-one, i.e. $g^*\in{\mathcal H}(0,m,n+1)$.  This completes the induction. 
\end{proof}

%%%%%%%%%%%%%%%%%%%%%%%%%%%%%%%%%

%%%%%%%%%%%%%%%%%%%%%%%%%%%%%%
%%%%%%%%%%%%%%%%%%%%%%%%%%%%%%
\vspace{0.5cm}

\section{Proof of Theorem 1.1 -  the general case}

By Corollary 2.5, we can assume that $p\leq m+k+1$. Representing $X$ as the union of an increasing sequence of closed sets  $X_i\subset X$ such that  each $f|X_i$ is perfect
and using that all restriction maps $\pi_i\colon C(X,\uin^{k+p})\to C(X_i,\uin^{k+p})$ are open and surjective, we can show that the proof of Theorem 1.1 is reduced to the case $f$ is a perfect map (see the proof of Proposition 2.5 for a similar situation). 
So, everywhere below we can suppose that the map $f$ from Theorem 1.1 is perfect.

Another reduction of Theorem 1.1 is provided by the following observation. By Lemma 2.4, the set 
${\mathcal H}_{\omega}(k,m,p)=C_{(m+k-p+2,\omega)}(X,\uin^{k+p})$ is open in $C(X,\uin^{k+p})$ for every 
$\omega\in cov(X)$. Since ${\mathcal H}(k,m,p)=\cap_{i=1}^{\infty}{\mathcal H}_{\omega_i}(k,m,p)$, where
$\{\omega_i\}\subset cov(X)$ is a sequence with $mesh(\omega_i)<2^{-i}$,
it suffices to show that  ${\mathcal H}_{\omega}(k,m,p)$ is dense in $C(X,\uin^{k+p})$ for every $\omega\in cov(X)$.
The remaining part of this section is devoted to the proof of this fact. We need few lemmas, in all these lemmas we suppose that $X$, $Y$, $f$ and the numbers $m,k,p$ are as in Theorem 1.1 with $f$ perfect. We also fix $\omega\in cov(X)$. 

\begin{lem}
If $C(X,\uin^{k+p})$ is equipped with the uniform convergence topology, then the set-valued map $\psi$ from $Y$ into 
$C(X,\uin^{k+p})$, defined by the formula
$\psi (y)=C(X,\uin^{k+p})\backslash C_{(m+k-p+2,\omega)}(X|f^{-1}(y),\uin^{k+p})$, has a closed graph.
\end{lem}

\begin{proof}
We can prove this lemma by following the arguments from the proof of  \cite[Lemma 2.6]{tv:02},
but in the present situation there exists a shorter proof.

Let $G=\cup\{y\times\psi (y):y\in Y\}\subset Y\times C(X,\uin^{k+p})$ be the graph of $\psi$ and $\{(y_n,g_n)\}$ a sequence in $G$ converging to $(y_0,g_0)\in Y\times C(X,\uin^{k+p})$.  It suffices to show that $(y_0,g_0)\in G$. 
Assuming $(y_0,g_0)\not\in G$, we conclude that $g_0\not\in\psi (y_0)$, so  
$g_0\in C_{(m+k-p+2,\omega)}(X|f^{-1}(y_0),\uin^{k+p})$. Then, by Lemma 2.2, there exists a
neighborhood $U$ of $y_0$ in $Y$ with $g_0|f^{-1}(U)$ being an $(m+k-p+2,\omega)$-map. We can suppose that $f^{-1}(y_n)\subset f^{-1}(U)$ for every $n$ because $\lim y_n=y_0$.  Consequently, $g_0|K$ is also an $(m+k-p+2,\omega)$-map, where $K$ denotes the union of all $f^{-1}(y_n)$,  $n=0, 1, 2,..$.  Obviously, $K$ is compact and, according to Lemma 2.4 (applied to the constant map  $q\colon K\to\{0\}$), the set $W$ of all $(m+k-p+2,\omega)$-maps $h\in C(K,\uin^{k+p})$ is open in $C(K,\uin^{k+p})$.  Since the sequence $\{g_n|K\}$ converges to $g_0|K$ in $C(K,\uin^{k+p})$ and $g_0|K\in W$, $g_n|K\in W$ for allmost all $n$. 
Therefore, there exists $j$ such that $g_j|f^{-1}(y_j)$ is an $(m+k-p+2,\omega)$-map.  The last conclusion contradicts the observation that $(y_j,g_j)\in G$ implies $g_j\not\in C_{(m+k-p+2,\omega)}(X|f^{-1}(y_j),\uin^{k+p})$.
Thus, $(y_0,g_0)\in G$. 
\end{proof}

Recall that a closed subset $F$ of the metrizable apace $M$ is said to be a $Z_n$-set in $M$,  where $n$ is a positive integer or 0,  if the set $C(\uin^n,M\backslash F)$ is dense in $C(\uin^n,M)$ with respect to the uniform convergence topology.

\begin{lem}
Let $\alpha\colon X\to (0,\infty)$ be a positive continuous function and $g_0\in C(X,\uin^{k+p})$.
Then $\psi (y)\cap\overline{B}(g_0,\alpha)$ is a $Z_m$-set in $\overline{B}(g_0,\alpha)$ for every $y\in Y$, where $\overline{B}(g_0,\alpha)$ is considered as a subspace of $C(X,\uin^{k+p})$ with the uniform convergence topology. 
\end{lem}

\begin{proof}
The proof of this lemma follows the proof of \cite[Lemma 2.8]{tv:02}. For sake of completeness we provide a sketch.
In this proof all function spaces are equipped with the uniform convergence topology generated by the Euclidean metric $d$ on $\uin^{k+p}$.
Since, by Lemma 3.1, $\psi$ has a closed graph, each $\psi (y)\cap\overline{B}(g_0,\alpha)$
is closed in $\overline{B}(g_0,\alpha)$. We need to show that, for fixed $y\in Y$, $\delta>0$ and a map $u\colon\uin^m \to \overline{B}(g_0,\alpha)$ there exists a map
$v\colon\uin^m\to\overline{B}(g_0,\alpha)\backslash\psi (y)$ which is $\delta$-close to $u$. Observe that $u$ generates $h\in C(\uin^m\times X,\uin^{k+p})$, $h(z,x)=u(z)(x)$, such that
$d(h(z,x),g_0(x))\leq\alpha (x)$ for any $(z,x)\in\uin^m\times X$. Since $f^{-1}(y)$ is compact, take $\lambda\in (0,1)$ such that $\lambda\sup\{\alpha (x):x\in f^{-1}(y)\}<\displaystyle\frac{\delta}{2}$ and define $h_1\in C(\uin^m\times f^{-1}(y),\uin^{k+p})$ by $h_1(z,x)=(1-\lambda)h(z,x)+\lambda g_0(x)$. Then, for every $(z,x)\in\uin^m\times f^{-1}(y)$, we have \\

\smallskip\noindent   
(1) \hbox{}~~~~~~$d(h_1(z,x),g_0(x))\leq (1-\lambda)\alpha (x)<\alpha (x)$ \\

\smallskip\noindent
and

\smallskip\noindent
(2) \hbox{}~~~~~~$d(h_1(z,x),h(z,x))\leq\lambda\alpha (x)<\displaystyle\frac{\delta}{2}$.

\smallskip\noindent
Let $\displaystyle q<\min\{r,\frac{\delta}{2}\}$, where $r=\inf\{\alpha (x)-d(h_1(z,x),g_0(x)):(z,x)\in\uin^m\times f^{-1}(y)\}$. 
Since $\dim f^{-1}(y)\leq k$, by Proposition 2.5 (applied to the projection $pr\colon\uin^m\times f^{-1}(y)\to\uin^m$),   
there is a map $h_2\in C(\uin^m\times f^{-1}(y),\uin^{k+p})$ such that $d(h_2(z,x),h_1(z,x))<q$ and $h_2|(\{z\}\times f^{-1}(y))$ is an
$(m+k-p+2,\omega)$-map for each $(z,x)\in\uin^m\times f^{-1}(y)$. Then, by $(1)$ and $(2)$, for all $(z,x)\in\uin^m\times f^{-1}(y)$ we have \\

\smallskip\noindent
(3) \hbox{}~~~~~~$d(h_2(z,x),h(z,x))<\delta$ and $d(h_2(z,x),g_0(x))<\alpha (x)$. \\

\smallskip\noindent
The equality $u_2(z)(x)=h_2(z,x)$ defines 
the map $u_2\colon\uin^m\to C(f^{-1}(y),\uin^{k+p})$. 
As in the proof of \cite[Lemma 2.8]{tv:02}, we can show that the map $\pi\colon\overline{B}(g_0,\alpha)\to C(f^{-1}(y),\uin^{k+p})$, $\pi (g)=g|f^{-1}(y)$, is continuous and open and
$u_2(z)\in\pi (\overline{B}(g_0,\alpha))$ for every $z\in\uin^m$.  So,
$\theta(z)=\overline{\pi^{-1}(u_2(z))\cap B_{\delta}(u(z))}$ 
 defines a convex-valued map from $\uin^m$ into $\overline{B}(g_0,\alpha)$ which is lower semi-continuous.
 Here, $B_{\delta}(u(z))$ is the open ball in $C(X,\uin^{k+p})$ (equipped with the uniform metric ) having center  $u(z)$  and radius $\delta$.  By  the Michael selection theorem \cite[Theorem 3.2"]{em:56}, there is a continuous selection $v\colon\uin^m\to C(X,\uin^{k+p})$ for $\theta$. Then $v$ maps $\uin^m$ into $\overline{B}(g_0,\alpha)$ and $v$ is $\delta$-close to $u$. Moreover, for any $z\in\uin^m$ we have $\pi(v(z))=u_2(z)$ and $u_2(z)$, being the restriction
$h_2|(\{z\}\times f^{-1}(y))$, is an $(m+k-p+2,\omega)$-map. Hence, 
$v(z)\not\in\psi (y)$ for any $z\in\uin^m$, i.e.
$v\colon\uin^m\to\overline{B}(g_0,\alpha)\backslash\psi (y)$.       
\end{proof}

Next lemma will finally accomplish the proof of Theorem 1.1.

\begin{lem}
The set ${\mathcal H}_{\omega}(k,m,p)$ is dense in $C(X,\uin^{k+p})$.
\end{lem}

\begin{proof}
Recall that by ${\mathcal H}_{\omega}(k,m,p)$ we denoted the set $C_{(m+k-p+2,\omega)}(X,\uin^{k+p})$.
It suffices to show that, for fixed $g_0\in C(X,\uin^{k+p})$ and a positive continuous function $\alpha\colon X\to (0,\infty)$, there exists $g\in \overline{B}(g_0,\alpha)\cap C_{(m+k-p+2,\omega)}(X,\uin^{k+p})$. To this end, consider the space $C(X,\uin^{k+p})$ with the uniform convergence topology as a closed and convex subset of the Banach space $E$ consisting of all bounded  maps from $X$ into ${\mathbb R}^{k+p}$. We 
define
the constant set-valued (and hence, lower semi-continuous) map $\phi$ from $Y$ into $C(X,\uin^{k+p})$, 
$\phi(y)=\overline{B}(g_0,\alpha)$. According to Lemma 3.2, 
$\overline{B}(g_0,\alpha)\cap\psi (y)$ is a $Z_m$-set in $\overline{B}(g_0,\alpha)$ for every $y\in Y$. So, we have a lower semi-continuous closed and convex-valued map $\phi$ from $Y$ to $E$ and a map $\psi\colon Y\to 2^E$ such that $\psi$ has a closed graph (see Lemma 3.1) and $\phi(y)\cap\psi (y)$ is a $Z_m$-set in $\phi(y)$ for each $y\in Y$.  Moreover, $\dim Y\leq m$, so we can apply \cite[Theorem 1.2]{gv:01} to obtain a  continuous map $h\colon Y\to E$ with $h(y)\in\phi(y)\backslash\psi (y)$ for every $y\in Y$.
Observe that $h$ is a map from $Y$ into $\overline{B}(g_0,\alpha)$ such that $h(y)\not\in\psi (y)$ for every $y\in Y$, i.e. 
$h(y)\in \overline{B}(g_0,\alpha)\cap C_{(m+k-p+2,\omega)}(X|f^{-1}(y),\uin^{k+p})$, $y\in Y$. Then
$g(x)=h(f(x))(x)$, $x\in X$, defines a map $g\in \overline{B}(g_0,\alpha)$ such that 
$g\in C_{(m+k-p+2,\omega)}(X|f^{-1}(y),\uin^{k+p})$ for every $y\in Y$. Hence, by virtue of Corollary 2.3, $g\in C_{(m+k-p+2,\omega)}(X,\uin^{k+p})$.  
\end{proof}

%%%%%%%%%%%%%%%%%%%%%%%%%%%%%%%%%

%%%%%%%%%%%%%%%%%%%%%%%%%%%%%%
%%%%%%%%%%%%%%%%%%%%%%%%%%%%%%

\bigskip

\end{document}